\documentclass[12pt]{article}
\usepackage{amssymb, amsmath}
\begin{document}
\renewcommand{\theequation}{\arabic{section}.\arabic{equation}}
\newtheorem{thm}{Theorem}[section]
\newtheorem{lem}{Lemma}[section]
\newtheorem{deff}{Definition}[section]
\newtheorem{pro}{Proposition}[section]
\newtheorem{cor}{Corollary}[section]
\newcommand{\n}{\nonumber}
\newcommand{\tv}{\tilde{v}}
\newcommand{\nV}{{\nabla V}}
\newcommand{\nv}{{\nabla v}}
\newcommand{\tw}{\tilde{\omega}}
\renewcommand{\wp}{W^{m,p}}
\newcommand{\kkp}{\frac{p+3}{kp}}
\newcommand{\kkpg}{\frac{kp\gamma}{p+3}}
\newcommand{\ppn}{\|_{L^p} ^{\frac{p+3}{kp}}}
\newcommand{\ppnn}{\|_{L^p} ^{1-\frac{p+3}{kp}}}
\newcommand{\kp}{\frac{p+2}{kp}}
\newcommand{\kpg}{\frac{kp\gamma}{p+2}}
\newcommand{\pn}{\|_{L^p} ^{\frac{p+2}{kp}}}
\newcommand{\pnn}{\|_{L^p} ^{1-\frac{p+2}{kp}}}
\renewcommand{\t}{\theta}
\newcommand{\T}{\Theta}
\newcommand{\nO}{\nabla \Theta}
\newcommand{\nOo}{\nabla^\bot \Theta}
\newcommand{\nt}{\nabla\theta}
\newcommand{\bn}{\|_{\dot{B}^0_{\infty, 1}}}
\newcommand{\w}{\omega}
\newcommand{\npg}{\frac{(p-3)\gamma}{2p}}
\newcommand{\np}{\frac{2p}{p-3}}
\newcommand{\npi}{\frac{p-3}{2p}}
\newcommand{\npn}{\|_{L^p}^{\frac{2p}{p-3}}}
\newcommand{\g}{\gamma}
\renewcommand{\k}{\kappa}
\newcommand{\e}{\varepsilon}
\renewcommand{\a}{\alpha}
\newcommand{\la}{\Lambda}
\newcommand{\laa}{\lambda}
\newcommand{\vare}{\varepsilon}
\newcommand{\s}{\sigma}
\newcommand{\hn}{\|_{\dot{H}^3}}
\renewcommand{\o}{\omega}
\renewcommand{\O}{\Omega}
\newcommand{\bb}{\begin{equation}}
\newcommand{\ee}{\end{equation}}
\newcommand{\bq}{\begin{eqnarray}}
\newcommand{\eq}{\end{eqnarray}}
\newcommand{\bqn}{\begin{eqnarray*}}
\newcommand{\eqn}{\end{eqnarray*}}
\title{On the  \emph{a priori} estimates for the Euler, the Navier-Stokes and the quasi-geostrophic equations}
\author{Dongho Chae\thanks{This research was supported partially
by  KRF Grant(MOEHRD, Basic Research Promotion Fund); $(\dagger)$ permanent address.\newline
}\\
Department of Mathematics\\
University of Chicago\\
Chicago, IL 60637, USA\\
  e-mail: {\it chae@math.uchicago.edu}\\
and \\
 Sungkyunkwan University$^{(\dagger)}$\\
Suwon 440-746, Korea}
 \date{}
\maketitle
\begin{abstract}
We prove new \emph{a priori}  estimates for the 3D Euler, the 3D Navier-Stokes and the 2D quasi-geostrophic equations by the method of similarity transforms.
\end{abstract}
\section{Main Results}
 \setcounter{equation}{0}

We are  concerned  on the following Navier-Stokes equations(Euler equations for $\nu=0$) describing the
homogeneous incompressible fluid flows in $\Bbb R^3$.
 \[
\mathrm{ (NS)_{\nu}}
 \left\{ \aligned
 &\frac{\partial v}{\partial t} +(v\cdot \nabla )v =-\nabla \textbf{p } +\nu \Delta v,
 \quad (x,t)\in {\Bbb R^3}\times (0, \infty) \\
 &\quad \textrm{div }\, v =0 , \quad (x,t)\in {\Bbb R^3}\times (0,
 \infty)\\
  &v(x,0)=v_0 (x), \quad x\in \Bbb R^3
  \endaligned
  \right.
  \]
where $v=(v_1, v_2, v_3 )$, $v_j =v_j (x, t)$, $j=1,2,3$, is the
velocity of the flow, $\mathbf{p}=\mathbf{p}(x,t)$ is the scalar pressure, $\nu \geq 0$ is the viscosity, and $v_0
$ is the given initial velocity, satisfying div $v_0 =0$. Given $m\in \Bbb N$,  we use  $W^{m,p}(\Bbb R^n )$ to denote the standard Sobolev space with the norm
$$\|f\|_{W^{m,p}}=\left(\sum_{|\a|\leq m }\int_{\Bbb R^n} |D^\a f(x)|^pdx \right)^{\frac1p},$$
where $\a=(\a_1, \cdots ,\a_n )$, $|\a|=\a_1 +\cdots +\a_n$ are the standard multi-index notation.
We also use $\dot{W}^{m,p} (\Bbb R^n )$ to denote the homogeneous space with the norm,
$$\|f\|_{\dot{W}^{m,p}}=\left(\sum_{|\a|= m }\int_{\Bbb R^n} |D^\a f(x)|^pdx \right)^{\frac1p}.$$
In the Hilbert space cases we denote $W^{m, 2}(\Bbb R^n)=H^m (\Bbb R^n)$, and $\dot{W}^{m, 2}(\Bbb R^n)=\dot{H}^m (\Bbb R^n)$.
The local well-posedness of the system $(NS)_\nu$ in $W^{m,p}(\Bbb R^3)$, $m>\frac{3}{p}+1$, is established in \cite{kat2, kat3}.
 The finite time blow-up problem (or equivalently the  regularity problem) of the local classical solution for both of the Euler equations and the Navier-Stokes are known as one of the most
 important and difficult problems in partial differential equations(see e.g. \cite{ler,caf} for the pioneering work and a later major advancement on the Navier-Stokes equations. see also \cite{maj,che,con1,cha1,lad,tem} for graduate level texts and survey articles on the current status of the problems for both of the Euler and the Navier-Stokes equations). The celebrated Beale-Kato-Majda criterion(\cite{bea}) states that the blow-up(for both of the Navier-Stokes and the Euler equations) happens at $T<\infty$ if and only if
$ \int_0 ^{T} \|\o (t )\|_{L^\infty} dt=\infty, $ where $\o =$ curl $v$ is the vorticity.
Motivated by Leray's question on the possibility of self-similar singularity in the Navier-Stokes equations(\cite{ler}), there are some nonexistence results  on the self-similar singularities for the Navier-Stokes equations(\cite{nec,tsa,mil}) and for the Euler equations(\cite{cha2,cha3,cha0}).
Transforming  the original Navier-Stokes and the Euler equations to the self-similar one(called the Leray equations in the case of Navier-Stokes equations), using appropriate similarity variables,  they  made analysis  the new system of equations to reach such nonexistence results.
In a recent preprint \cite{cha0}, new type of similarity transforms  which depend on the solution itself  are considered, and with suitable choice of its form  some of \emph{a priori}  estimates  are derived for the smooth solutions of the Euler and the Navier-Stokes equations. The purpose of this paper is develop further the method to prove  high order derivative estimates for the Euler,  the Navier-Stokes equations and also for the quasi-geostrophic equations as well as the general $L^p$ estimates for the Navier-Stokes equations. In the quasi-geostrophic equations for the critical space case we need to use critical Besov spaces, and the  derivation of estimates rely on the particle trajectory method for the transformed system. We state our main theorems below.

\begin{thm} Suppose $m\geq 3  $, $m\in \Bbb N$, be given.
Let $v_0 \in H^m (\Bbb R^3)$, and $v\in C([0, T); H^m (\Bbb R^3))$ be the classical solution of the
system $(NS)_\nu$. Then,  for all $t\in [0, T)$ and $k\in \{3, \cdots, m\}$ there exists $C_0 =C_0 (k)$ such that for all $\g \geq C_0  \|v_0\|_{L^2}^{1-\frac{5}{2k}}$ the following inequalities holds true:\\
\ \\
\noindent(i) for  the case $\nu\geq 0$, we have
\bb\label{thm1}
\|D^k v(t)\|_{L^2} \leq \frac{\|D^k v_0\|_{L^2} \exp \left[ \frac{2k\g}{5} \int_0 ^t \|D^k v (\tau)\|_{L^2} ^{\frac{5}{2k}} d\tau\right]}{\left\{1+\left(\g-C_0   \|v_0\|_{L^2}^{1-\frac{5}{2k}}\right)
\|D^k v_0 \|_{L^2}^{\frac{5}{2k}}\int_0 ^t \exp\left[\g \int_0 ^\tau \|D^k v(\sigma )\|_{L^2}^{\frac{5}{2k}}d\sigma\right] d\tau \right\}^{\frac{2k}{5}}},
\ee
with an upper estimate of the denominator,
\bq\label{thm1a}
\lefteqn{1+\left(\g-C_0  \|v_0\|_{L^2}^{1-\frac{5}{2k}}\right)
\|D^k v_0 \|_{L^2}^{\frac{5}{2k}}\int_0 ^t \exp\left[\g \int_0 ^\tau \|D^k v(\sigma )\|_{L^2}^{\frac{5}{2k}}d\sigma\right] d\tau}\hspace{1.in}\n \\
&&\leq \frac{1}{\left(1- C_0 \|v_0\|_{L^2}^{1-\frac{5}{2k}}\|D^k v_0 \|_{L^2}^{\frac{5}{2k}} t\right)^{\frac{\g}{C_0\|v_0\|_{L^2}^{1-\frac{5}{2k}}}-1}}.\n \\
\eq

\noindent(ii)  for the case  $\nu=0$, we have
\bb\label{thm2}
\|D^k v(t)\|_{L^2} \geq \frac{\|D^k v_0\|_{L^2} \exp \left[- \frac{2k\g}{5} \int_0 ^t \|D^k v (\tau)\|_{L^2} ^{\frac{5}{2k}} d\tau\right]}{\left\{1-\left(\g-C_0   \|v_0\|_{L^2}^{1-\frac{5}{2k}}\right)
\|D^k v_0 \|_{L^2}^{\frac{5}{2k}}\int_0 ^t \exp\left[-\g \int_0 ^\tau \|D^k v(\sigma )\|_{L^2}^{\frac{5}{2k}}d\sigma\right] d\tau \right\}^{\frac{2k}{5}}}
\ee
with a lower estimate of the denominator,
\bq\label{thm2a}
\lefteqn{1-\left(\g-C_0  \|v_0\|_{L^2}^{1-\frac{5}{2k}}\right)
\|D^k v_0 \|_{L^2}^{\frac{5}{2k}}\int_0 ^t \exp\left[-\g \int_0 ^\tau \|D^k v(\sigma )\|_{L^2}^{\frac{5}{2k}}d\sigma\right] d\tau}\hspace{1.in}\n \\
&&\geq \frac{1}{\left(1+ C_0 \|v_0\|_{L^2}^{1-\frac{5}{2k}}\|D^k v_0 \|_{L^2}^{\frac{5}{2k}} t\right)^{\frac{\g}{C_0\|v_0\|_{L^2}^{1-\frac{5}{2k}}}-1}}.\n \\
\eq
\end{thm}

\noindent{\it Remark 1.1 } In the special case $\g =C_0 \|v_0\|_{L^2}^{1-\frac{5}{2k}}$ the estimates (\ref{thm1}) and (\ref{thm2}) reduces to the form, which could be also be proved directly from $(NS)_\nu$ without using the similarity transform as in the proof below. The main novelty of the above estimates and all the other estimates in the theorems below is that  $\g$ is a free parameter that can take any value greater or equal to a constant, which makes nontrivial increment in time of the denominator in (\ref{thm1})(decrement of the denominator in (\ref{thm2})).  An interesting problem to consider is  `optimization' of those estimates by suitable choice of $\g$.\\
\ \\
\noindent{\it Remark 1.2 } The estimate (\ref{thm2a}) shows that the finite time blow-up of the Euler equations, even if it is true, does not follow from the inequality (\ref{thm2}).\\
\ \\
In the following theorem we restrict $\nu >0$, hence it is only for the Navier-Stokes equations.
Before its statement we recall that the local in time well-posedness in $L^p(\Bbb R^3)$ of the Navier-Stokes equations is proved by \cite{kat1}.

\begin{thm} Let  $p\in (3, \infty)$ be given.
Suppose  $v_0 \in L^p(\Bbb R^3)$, and $v\in C([0, T); L^p (\Bbb R^3))$ be the classical solution of the
system $(NS)_\nu$, $\nu >0$. Then,  for all $t\in [0, T)$ there exists $C_0 =C_0 (\nu, p)$ such that for all $\g \geq C_0 $ the following inequality holds true:\\
\bb\label{thm7}
\| v(t)\|_{L^p} \leq \frac{\| v_0\|_{L^p} \exp \left[ \frac{(p-3)\g}{2p}  \int_0 ^t \| v (\tau)\|_{L^p} ^{\frac{2p}{p-3}} d\tau\right]}{\left\{1+\left(\g-C_0 \right)
\| v_0 \|_{L^p}^{\frac{2p}{p-3}}\int_0 ^t \exp\left[\g \int_0 ^\tau \| v(\sigma )\|_{L^p}^{\frac{2p}{p-3}}d\sigma\right] d\tau \right\}^{\frac{p-3}{2p}}},
\ee
with an upper estimate of the denominator,
\bq\label{thm7a}
\lefteqn{1+\left(\g-C_0\right)
\| v_0 \|_{L^p}^{\frac{2p}{p-3}}\int_0 ^t \exp\left[\g \int_0 ^\tau \|v(\sigma )\|_{L^p}^{\frac{2p}{p-3}}d\sigma\right] d\tau}\hspace{1.in}\n \\
&&\leq \frac{1}{\left(1- C_0\|v_0 \|_{L^p}^{\frac{2p}{p-3}} t\right)^{\frac{\g}{C_0}-1}}.\n \\
\eq
\end{thm}

Next we are concerned on deriving  estimates for the  two dimensional dissipative
quasi-geostrophic equation:
\[
(QG)_{\kappa}\left\{\aligned
 &\partial_t \theta +(v\cdot \nabla ) \theta+\kappa \la ^{\alpha}
 \theta=0,\quad
 (x,t)\in {\Bbb R^3}\times (0,
 \infty),\quad\alpha \geq 0,\\
 & v = \nabla^{\perp} (-\Delta)^{-\frac12}\theta, \\
&\theta(0,x)=\theta_0,
\endaligned \right.\]
where $\la =(-\Delta )^{\frac12}.$
After pioneering work by Constantin-Majda-Tabak(\cite{con4} the system $(QG)_{\kappa}$ became a hot subject of studies(see e.g. \cite{con3a, wu1, cor0, caf1, kis} and references therein), mainly due to its structural resemblance to the 3D Euler and the 3D Navier-Stokes equations with similar difficulties in the regularity problems.
Contrary to the case of the system $(NS)_\nu$, where we only have control of $L^2$ norm of velocity,
we have the following $L^p$ bound of $\t$ for any $p\in [1, \infty]$ in the system $(QG)_\k$,
$$ \|\theta (t)\|_{L^p}\leq \|\theta _0\|_{L^p}. $$
Due to this fact we can apply our method to derive $W^{k,p}$  estimates for $(QG)_\k$ as follows.
\begin{thm} Suppose $1<p<\infty$,  and $m > \frac{2}{p}+1$, $m\in \Bbb N$, be given.
Let $\theta_0 \in W^{m,p} (\Bbb R^2)$, and $\theta \in C([0, T); W^{m,p} (\Bbb R^2))$ be the classical solution of the
system $(QG)_\kappa$. Then, for all $t\in [0, T)$, $p\in (1, \infty)$, $k\in \{ [\frac{2}{p}+1] , \cdots, m\}$  there exists $C_0 =C_0 (k,p)$ such that and $\g \geq C_0  \|\theta _0\|_{L^p}^{1-\frac{p+2}{kp}}$ the following inequalities holds true:\\
\ \\
\noindent(i) for the case $\kappa \geq 0$ and $\a \geq 0$,
\bb\label{thm3}
\|D^k \t(t)\|_{L^p} \leq \frac{\|D^k \t_0\|_{L^p} \exp \left[ \frac{kp\g}{p+2} \int_0 ^t \|D^k \t (\tau)\|_{L^p} ^{\frac{p+2}{kp}} d\tau\right]}{\left\{1+\left(\g-C_0   \|\t_0\|_{L^p}^{1-\frac{p+2}{kp}}\right)
\|D^k \t_0 \|_{L^p}^{\frac{p+2}{kp}}\int_0 ^t \exp\left[\g \int_0 ^\tau \|D^k \t(\sigma )\|_{L^p}^{\frac{p+2}{kp}}d\sigma\right] d\tau \right\}^{\frac{kp}{p+2}}}
\ee
with an upper estimate of the denominator,
\bq\label{thm3a}
\lefteqn{1+\left(\g-C_0  \|\theta _0\|_{L^p}^{1-\frac{p+2}{kp}}\right)
\|D^k \t_0 \|_{L^p}^{1-\frac{p+2}{kp}}\int_0 ^t \exp\left[\g \int_0 ^\tau \|D^k \t (\sigma )\|_{L^p}^{\frac{p+2}{kp}}d\sigma\right] d\tau}\hspace{1.in}\n \\
&&\leq \frac{1}{\left(1- C_0 \|\t_0\|_{L^p}^{1-\frac{p+2}{kp}}\|D^k \t_0 \|_{L^p}^{\frac{p+2}{kp}} t\right)^{\frac{\g}{C_0\|\t_0\|_{L^p}^{1-\frac{p+2}{kp}}}-1}}.\n \\
\eq

\noindent(ii) for the case  $\kappa=0$,
\bb\label{thm4}
\|D^k \t(t)\|_{L^p} \geq\frac{\|D^k \t_0\|_{L^p} \exp \left[ -\frac{kp\g}{p+2} \int_0 ^t \|D^k \t (\tau)\|_{L^p} ^{\frac{p+2}{kp}} d\tau\right]}{\left\{1-\left(\g-C_0   \|\t_0\|_{L^p}^{1-\frac{p+2}{kp}}\right)
\|D^k \t_0 \|_{L^p}^{\frac{p+2}{kp}}\int_0 ^t \exp\left[-\g \int_0 ^\tau \|D^k \t(\sigma )\|_{L^p}^{\frac{p+2}{kp}}d\sigma\right] d\tau \right\}^{\frac{kp}{p+2}}}
\ee
with a lower estimate of the denominator,
\bq\label{thm4a}
\lefteqn{1-\left(\g-C_0  \|\theta _0\|_{L^p}^{1-\frac{p+2}{kp}}\right)
\|D^k \t_0 \|_{L^p}^{1-\frac{p+2}{kp}}\int_0 ^t \exp\left[-\g \int_0 ^\tau \|D^k \t (\sigma )\|_{L^p}^{\frac{p+2}{kp}}d\sigma\right] d\tau}\hspace{1.in}\n \\
&&\geq \frac{1}{\left(1+C_0 \|\t_0\|_{L^p}^{1-\frac{p+2}{kp}}\|D^k \t_0 \|_{L^p}^{\frac{p+2}{kp}} t\right)^{\frac{\g}{C_0\|\t_0\|_{L^p}^{1-\frac{p+2}{kp}}}-1}}.\n \\
\eq
\end{thm}

In the critical space case with $m \simeq\frac{2}{p}+1$,  we have different form of estimate,
where  the  use of  the critical Besov space, $\dot{B}^0_{\infty,1}$, is necessary. For a brief introduction of this Besov space please see the next section.

\begin{thm}
Let $\theta \in C([0, T); \dot{B}^1_{\infty,1} (\Bbb R^2))$ be a classical solution of $(QG)_0$
with initial data $\t_0 \in \dot{B}^1_{\infty,1} (\Bbb R^2)$, then there exists $C_0$ such that for all $t\in [0, T)$ and $\g \geq C_0$ we have the following  upper estimate,
\bq\label{thm5}
\|\nt (t)\|_{L^\infty}\leq \frac{\|\nt_0\|_{L^\infty} \exp\left[ \g \int_0 ^t \|\nt (\tau)\bn d\tau\right]}{1+ (\g-C_0)\|\nt_0 \|_{L^\infty}\int_0 ^t\exp\left[  \g\int_0 ^\tau\|\nt(\sigma)\bn d\sigma \right]d\tau},\n\\
 \eq
 and lower one
  \bq\label{thm6}
\|\nt (t)\|_{L^\infty}\geq \frac{\|\nt_0\|_{L^\infty} \exp\left[ -\g \int_0 ^t \|\nt (\tau)\bn d\tau\right]}{1-(\g-C_0) \|\nt_0 \|_{L^\infty}\int_0 ^t\exp\left[ -\g \int_0 ^\tau\|\nt(\sigma)\bn d\sigma \right]d\tau}.\n \\
 \eq
In particular, the denominator of the right hand side of  (\ref{thm6}) can be  estimated from below as follows.
  \bq\label{thm6a}
\lefteqn{1-(\g-C_0 )
\|\nt_0 \|_{L^\infty}\int_0 ^t \exp\left[-\g  \int_0 ^\tau \|\nt (\sigma )\bn d\sigma\right] d\tau}\hspace{2.in}\n \\
&&\geq \frac{1}{\left(1+C_0 \|\nt_0\|_{L^\infty} t\right)^{\frac{\g}{C_0}-1}}.
\eq
\end{thm}
{\it Remark 1.3 } As will be seen clearly  in the proof below, the optimal constant $C_0$ in the above theorem is the optimal constant in the following Calderon-Zygmund type of inequality,
$$
\|\nabla v\bn \leq C_0 \|\nt \bn.
$$
\noindent{\it Remark 1.4 } In the special case of $\g=C_0$ the above estimates (\ref{thm5}) and (\ref{thm6})  reduce to the well known ones that could be directly obtained from $(QG)_0$ by the standard method.

\section{Proof of the Main Results}
\setcounter{equation}{0}

We first recall the following well-known inequalities:
\begin{itemize}
\item[(a)] For $k>\frac{n}{p}+1$ and  $f,g\in W^{k,p} (\Bbb R^n)$ there exists constant $C_1=C_1 (k,p,n)$ such that
\bb\label{com}
\|D^k (fg)-fD^k g\|_{L^p}\leq C_1 \left(\|\nabla f\|_{L^\infty} \|D^{k-1} g\|_{L^p}
+\|D^k f\|_{L^p} \|g\|_{L^\infty}\right).
\ee
(the commutator estimate, \cite{kla, kat3})
\item[(b)] For $k>\frac{n}{p}+1$ and  $f,g\in W^{k,p} (\Bbb R^n)$ there exists constant $C_2=C_2 (k,p,n)$ such that
\bb\label{gag}
\|\nabla f\|_{L^\infty} \leq C_2 \|f\|_{L^p} ^{1-\frac{p+n}{kp}} \|D^k f\|_{L^p} ^{\frac{p+n}{kp}}.
\ee
(the Gagliardo-Nirenberg inequality, \cite{ada})
\end{itemize}
For $\a \in [0, 2]$ we also recall the following estimate for the fractional laplacian
\bb\label{frac}
\int_{\Bbb R^n} |f|^{p-2} f\la ^\a f dx \geq \frac{2}{p} \int_{\Bbb R^n}
 \left(\la ^{\frac{\a}{2}} |f|^{\frac{p}{2}} \right)^2 dx.
 \ee
 (see\cite{ju} for the proof, and  see also \cite{cor0} for its earlier version),
 Below we briefly introduce some of the critical Besov spaces, which is necessary for our purpose(see e.g. \cite{tri} for more comprehensive introduction).
 Given $f\in \mathcal{S}$, the Schwartz class of rapidly deceasing
 functions, its Fourier transform $\hat{f}$ is defined by
 $$
 \mathcal{F}(f)= \hat{f} (\xi)=\frac{1}{(2\pi )^{n/2}}\int_{\Bbb R^n} e^{-ix\cdot \xi }
 f(x)dx.
  $$
 We consider  $\varphi \in \mathcal{S}$ satisfying the following three conditions:
\begin{itemize}
 \item[(i)] $\textrm{Supp}\, \hat{\varphi} \subset
 \{\xi \in {\Bbb R}^n  \,| \,\frac12 \leq |\xi|\leq
 2\}$,
 \item[(ii)] $\hat{\varphi} (\xi)\geq C >0 $ if $\frac23 <|\xi|<\frac32$,
 \item[(iii)] $
 \sum_{j\in \mathbb{Z}}  \hat{\varphi}_j (\xi )=1$,
  where $\hat{\varphi_j } =\hat{\varphi } (2^{-j} \xi
 )$.
 \end{itemize}
Construction  of such sequence of functions $\{ \varphi_j \}_{j\in
\Bbb Z}$ is well-known.
For $s\in \Bbb R$, space $\dot{B}^s_{\infty ,
 1}$ is defined by
 $$ f\in \dot{B}^s_{\infty ,
 1}\Longleftrightarrow \|f\|_{\dot{B}^s_{\infty ,
 1}}=\sum_{j\in \Bbb Z}2^{sj}\|\varphi_j *
 f\|_{L^\infty} < \infty,
 $$
 where $*$ is the standard notation for convolution,
 $(f*g)(x)=\int_{\Bbb R^n} f(x-y)g(y)dy$. The norm $\|\cdot
 \|_{\dot{B}^s_{\infty,1}}$  is actually defined up to addition of
 polynomials(namely,
  if $f_1-f_2 $ is a polynomial, then both of $f_1$ and $f_2$ give the same norm),
 and the space  $\dot{B}^s_{\infty,1} (\Bbb R^n)$ is defined as the
 quotient space of  a class of functions with finite norm, $\|\cdot\|_{\dot{B}^s_{\infty,1}}$,
 divided by the space of polynomials in $\Bbb R^n$.
Note that the condition (iii) implies
immediately
\bq
 \|f\|_{L^\infty} \leq \|f\bn.
 \eq
The crucial feature of $\dot{B}^0_{\infty,1}(\Bbb R^n) $, compared with $L^\infty(\Bbb R^n)$ is that the singular
integral operators of the Calderon-Zygmund type map $\dot{B}^0_{\infty,1} (\Bbb R^n)$ into itself
boundedly, the property which $L^\infty$ does not have. See \cite{bou} for more details on these homogeneous Besov spaces.\\

\noindent{\bf Proof  of Theorem 1.1 }
 Let $T$ be the maximal time of existence of a classical solution $v$ of $(NS)_\nu$ in $\wp (\Bbb R^3)$, and $v\in C([0, T); \wp (\Bbb R^3 ))$.
Given a classical solution $v(x,t)$ and the associated pressure function $p(x,t)$, we
introduce a functional transform from
$(v,\textbf{p})$ to $(V, P)$ defined by the formula,
\bq\label{2.4}
 v(x,t)&=&\exp\left[\pm\frac{3\g}{5}\int_0 ^t\| D^k v(\tau)\|_{L^2} ^{\frac{5}{2k}} d\tau \right] \ V\left(y, s \right),\\
 \label{2.5}
 \textbf{p}(x,t)&=&\exp\left[\pm\frac{ 6\g}{5}\int_0 ^t\|D^k v(\tau)\|_{L^2}^{\frac{5}{2k}} d\tau \right] \ P\left(y, s\right)
 \eq
 with
\bq \label{2.6}
y&=&\exp\left[\pm \frac{2\g}{5}\int_0 ^t\|D^k v(\tau)\|_{L^2} ^{\frac{5}{2k}} d\tau \right]x,  \\
\label{2.6a}
 s&=&\int_0 ^t\exp\left[\pm \g\int_0 ^\tau\| D^k v(\sigma)\|_{L^2} ^{\frac{5}{2k}} d\sigma \right]d\tau ,
\eq
respectively for $(\pm)$. We note that this choice of similarity transform makes the scaling dimension of the energy become zero, and thus the energy invariant of the transform,
$$ \|v(t)\|_{L^2}=\|V(s)\|_{L^2}.
$$
We also note the following integral invariant of the transformation (\ref{2.4})-(\ref{2.6a}),
 $$
 \int_0 ^t \|D^k v(\tau )\|_{L^2}^{\frac{5}{2k}} d\tau= \int_0 ^s \|D^k V(\sigma)\|_{L^2}^{\frac{5}{2k}}d\sigma   \qquad 3\leq k <\infty.
 $$
Substituting $(v,\textbf{p})$ in (\ref{2.4})-(\ref{2.6a}) into  $(NS)_\nu$, we  obtain an equivalent system of equations:
$$
 (NS)_\nu ^\pm \left\{ \aligned
 &\mp\frac{\g }{5} \|D^k V(s)\|_{L^2} ^{\frac{5}{2k}} \left[ 3 V + 2 (y \cdot \nabla)V \right]=V_s+(V\cdot \nabla )V +\nabla P \\
& \qquad \hspace{1.in} -\nu \Delta V \exp \left[ \mp\frac{\g}{5}\int_0 ^s \|D^k V(\sigma)\|_{L^p}^{\frac{5}{2k}}d\sigma \right],\\
 & \mathrm{div}\, V=0,\\
&V(y,0)=V_0(y)=v_0 (y),
\endaligned \right.
$$
  where $(NS)_\nu ^+$ means that we have chosen (+) sign in (\ref{2.4})-(\ref{2.6a}), and this corresponds to $(-)$  sign in the first equations of $(NS)_\nu ^\pm $.  Similarly for $(NS)_\nu ^-$.

 We observe that $V\in C([0, S_\pm ); H^m (\Bbb R^3))$, where
 $$S_\pm:=\int_0 ^T\exp\left[\pm\g\int_0 ^\tau\| D^k v(\sigma)\|_{L^2} ^{\kkp} d\sigma \right]d\tau
 $$
 is the maximal time of existence of the classical solution in $H^m (\Bbb R^3 )$ for the system $(NS)_\nu ^\pm$ respectively.
 Form now on we separate our proof.\\

 \noindent{\it \underline{Proof of (i)}: } We choose $(+)$ sign in (\ref{2.4})-(\ref{2.6a}), and
 work with $(NS)_\nu ^+$, where $\nu \geq 0$.
 Taking $L^2 (\Bbb R^3)$ inner product of the first equations of $(NS)_\nu ^+$ with $V$, and integrating by part, we find that
 $$\frac12\frac{d}{ds} \|V(s)\|_{L^2}^2+\nu  \int_{\Bbb R^3} | \nabla V|^2 dy\,
 \exp \left[ -\frac{\g}{5}\int_0 ^s \|D^k V(\sigma)\|_{L^2}^{\frac{5}{2k}}d\sigma \right] =0.$$
Hence we have energy bound,
\bb\label{energy}
\|V(s)\|_{L^2}\leq \|V_0 \|_{L^2}.
\ee
Next, taking $\dot{H}^k (\Bbb R^3)$ inner product of  the first equations of $(NS)_\nu ^+ $ by $V$, and integrating by part,
we derive
\bq\label{est1}
\lefteqn{\frac12 \frac{d}{ds} \|D^k V\|_{L^2} ^2 +\frac{2k\g}{5} \|D^k V\|_{L^2} ^{2+\frac{5}{2k}} +\nu \|D^{k+1} V\|_{L^2}^2\exp \left[- \frac{\g}{5}\int_0 ^s \|D^k V(\sigma)\|_{L^2}^{\frac{5}{2k}}d\sigma \right]}\hspace{.9 in}\n \\
&&= - (D^k (V\cdot \nabla ) V- (V\cdot \nabla )D^k V, D^k V)_{L^2}\n \\
&& \leq C \|\nabla V\|_{L^\infty} \|D^k V\|_{L^2} ^2\leq C \|V\|_{L^2}^{1-\frac{5}{2k}} \|D^k V\|_{L^2} ^{2+\frac{5}{2k}}\n \\
&&\leq \frac{2kC_0}{5}\|V_0\|_{L^2}^{1-\frac{5}{2k}} \|D^k V\|_{L^2} ^{2+\frac{5}{2k}}
\eq
for an absolute constant $C_0=C_0 (k)$, where we used the  computations,
\bqn
\lefteqn{\left(D^k (y\cdot \nabla )V, D^k V\right)_{L^2}= \frac12 \int_{\Bbb R^3} (y\cdot \nabla )|D^k V|^2 dy + k\|D^k V\|_{L^2}^2}\hspace{.1in}\\
&&=-\frac32 \|D^k V\|_{L^2}^2 + k\|D^k V\|_{L^2}^2
=(k-\frac32) \|D^k V\|_{L^2}^2,
\eqn
the commutator estimate (\ref{com}) and the Gagliardo-Nirenberg inequality (\ref{gag}).
Hence, from (\ref{est1}), ignoring the viscosity term, we have the differential inequality
$$
\frac{d}{ds} \|D^k V\|_{L^2}\leq -\frac{2k}{5}\left(\g -C_0\|V_0\|_{L^2}^{1-\frac{5}{2k}}\right) \|D^k V\|_{L^2} ^{1+\frac{5}{2k}},
$$
which can be solved to provide us with
\bb\label{sol1}
\|D^k V(s)\|_{L^2}\leq \frac{\|D^k V_0\|_{L^2}}{\left[1+\left(\g-C_0 \|V_0\|_{L^2}^{1-\frac{5}{2k}}\right)
\|D^k V_0 \|_{L^2}^{\frac{5}{2k}} s\right]^{\frac{2k}{5}} }
\ee
for all $s\in [0, S_+)$.
Transforming back to the original velocity $v$, using the relations (\ref{2.4})-(\ref{2.6a}) with $(+)$ sign,  we  obtain (\ref{thm1}).
In order to derive (\ref{thm1a}) we observe that (\ref{thm1}) can be written in the integrable form,
\bq\label{log1}
&&\|D^k v(t)\|_{L^2}^{\frac{5}{2k}} \leq \frac{\|D^k v_0\|_{L^2}^{\frac{5}{2k}} \exp \left[ \g \int_0 ^t \|D^k v (\tau)\|_{L^2} ^{\frac{5}{2k}} d\tau\right]}{\left\{1+\left(\g-C_0  \|v_0\|_{L^2}^{1-\frac{5}{2k}}\right)
\|D^k v_0 \|_{L^2}^{\frac{5}{2k}}\int_0 ^t \exp\left[\g \int_0 ^\tau \|D^k v(\sigma )\|_{L^2}^{\frac{5}{2k}}d\sigma\right] d\tau \right\}}\n \\
&&=\left(\g -C_0  \|v_0\|_{L^2}^{1-\frac{5}{2k}}\right)^{-1} \times \n \\
&&\times \frac{d}{dt}\log \left\{1+\left(\g-C_0 \|v_0\|_{L^2}^{1-\frac{5}{2k}}\right)\|D^k v_0 \|_{L^2}^{\frac{5}{2k}}\int_0 ^t \exp\left[\g \int_0 ^\tau \|D^k v(\sigma )\|_{L^2}^{\frac{5}{2k}}d\sigma\right] d\tau \right\}.\n \\
\eq
Hence, integrating  (\ref{log1}) over $[0, t]$, we obtain
\bq\label{diff1}
&&\int_0 ^t\|D^k v(\tau)\|_{L^2}^{\frac{5}{2k}}d\tau\leq \left(\g-C_0 \|v_0\|_{L^2}^{1-\frac{5}{2k}}\right)^{-1}\times\n \\
&&\times \log \left\{1+\left(\g-C_0   \|v_0\|_{L^2}^{1-\frac{5}{2k}}\right)
\|D^k v_0 \|_{L^2}^{\frac{5}{2k}}\int_0 ^t \exp\left[\g \int_0 ^\tau \|D^k v(\sigma )\|_{L^2}^{\frac{5}{2k}}d\sigma\right] d\tau \right\}.\n \\
\eq
Now, setting
$$y(t):= 1+\left(\g-C_0  \|v_0\|_{L^2}^{1-\frac{5}{2k}}\right)
\|D^k v_0 \|_{L^2}^{\frac{5}{2k}}\int_0 ^t \exp\left[\g \int_0 ^\tau \|D^k v(\sigma )\|_{L^2}^{\frac56}d\sigma\right] d\tau,
$$
we find that (\ref{diff1}) can be rewritten as a differential  inequality,
\bb\label{diffeq1}
y'(t)\leq \left(\g-C_0  \|v_0\|_{L^2}^{1-\frac{5}{2k}}\right)
\|D^k v_0 \|_{L^2}^{\frac{5}{2k}} \, y(t)^{M},
\ee
where we set
\bb\label{m1}
M:=\frac{\g}{\g-C_0  \|v_0\|_{L^2}^{1-\frac{5}{2k}}}.
\ee
 The differential inequality  (\ref{diffeq1}) is solved as
\bb\label{y1}
y(t)\leq \frac{1}{\left(1- C_0 \|v_0\|_{L^2}^{1-\frac{5}{2k}}\|D^k v_0 \|_{L^2}^{\frac{5}{2k}} t\right)^{\frac{\g}{C_0\|v_0\|_{L^2}^{1-\frac{5}{2k}}}-1}},
\ee
which provides us with (\ref{thm1a}).\\
\ \\
\noindent{\it \underline{Proof of (ii)}: } Here we choose $(-)$ sign in (\ref{2.4})-(\ref{2.6a}), and
 work with $(NS)_0 ^-$.
 Taking $L^2 (\Bbb R^3)$ inner product of the first equations of $(NS)_0 ^-$ with $V$, and integrating by part, we find that
 $$
\frac{d}{ds} \|V(s)\|_{L^2}^2 =0
 $$
  which implies  energy equality
\bb
 \|V(s)\|_{L^2}\leq \|V_0\|_{L^2},
 \ee
 Next, taking $\dot{H}^k (\Bbb R^3)$ inner product of the first equations of $(NS)_0 ^-$ with $V$, and integrating by part,
we derive similarly to the above
\bq\label{est2}
\lefteqn{\frac12 \frac{d}{ds} \|D^k V\|_{L^2} ^2 -\frac{2k\g}{5} \|D^k V\|_{L^2} ^{2+\frac{5}{2k}} }\hspace{.9 in}\n \\
&&\geq - \frac{2kC_0 }{5} \|V_0\|_{L^2}^{1-\frac{5}{2k}} \|D^k V\|_{L^2} ^{2+\frac{5}{2k}}
\eq
for the same absolute constant $C_0=C_0 (k)$ as in (\ref{est1}). Hence,
$$
\frac{d}{ds} \|D^k V\|_{L^2}\geq \frac{2k}{5}\left(\g -C_0\|V_0\|_{L^2}^{1-\frac{5}{2k}}\right) \|D^k V\|_{L^2} ^{1+\frac{5}{2k}},
$$
which can be solved to provide us with
\bb\label{sol2}
\|D^k V(s)\|_{L^2}\geq \frac{\|D^k V_0\|_{L^2}}{\left[1-\left(\g-C_0 \|V_0\|_{L^2}^{1-\frac{5}{2k}}\right)
\|D^k V_0 \|_{L^2}^{\frac{5}{2k}} s\right]^{\frac{2k}{5}} }
\ee
for all $s\in [0, S_-)$.
Transforming back to the original velocity $v$, using the relations (\ref{2.4})-(\ref{2.6a}) with $(-)$ sign,  we  have (\ref{thm2}).
In order to derive (\ref{thm2a}) we  rewrite (\ref{thm2}) in the integrable form,
\bq\label{log2}
&&\|D^k v(t)\|_{L^2}^{\frac{5}{2k}} \geq \frac{\|D^k v_0\|_{L^2}^{\frac{5}{2k}} \exp \left[ -\g \int_0 ^t \|D^k v (\tau)\|_{L^2} ^{\frac{5}{2k}} d\tau\right]}{\left\{1-\left(\g-C_0  \|v_0\|_{L^2}^{1-\frac{5}{2k}}\right)
\|D^k v_0 \|_{L^2}^{\frac{5}{2k}}\int_0 ^t \exp\left[-\g \int_0 ^\tau \|D^k v(\sigma )\|_{L^2}^{\frac{5}{2k}}d\sigma\right] d\tau \right\}}\n \\
&&=-\left(\g -C_0  \|v_0\|_{L^2}^{1-\frac{5}{2k}}\right)^{-1} \times \n \\
&&\times \frac{d}{dt}\log \left\{1-\left(\g-C_0 \|v_0\|_{L^2}^{1-\frac{5}{2k}}\right)\|D^k v_0 \|_{L^2}^{\frac{5}{2k}}\int_0 ^t \exp\left[-\g \int_0 ^\tau \|D^k v(\sigma )\|_{L^2}^{\frac{5}{2k}}d\sigma\right] d\tau \right\}.\n \\
\eq
Integrating  (\ref{log2}) over $[0, t]$, we obtain
\bq\label{diff2}
&&\int_0 ^t\|D^k v(\tau)\|_{L^2}^{\frac{5}{2k}}d\tau\geq -\left(\g-C_0 \|v_0\|_{L^2}^{1-\frac{5}{2k}}\right)^{-1}\times\n \\
&&\times \log \left\{1-\left(\g-C_0   \|v_0\|_{L^2}^{1-\frac{5}{2k}}\right)
\|D^k v_0 \|_{L^2}^{\frac{5}{2k}}\int_0 ^t \exp\left[-\g \int_0 ^\tau \|D^k v(\sigma )\|_{L^2}^{\frac{5}{2k}}d\sigma\right] d\tau \right\}.\n \\
\eq
Setting
$$y(t):= 1-\left(\g-C_0  \|v_0\|_{L^2}^{1-\frac{5}{2k}}\right)
\|D^k v_0 \|_{L^2}^{\frac{5}{2k}}\int_0 ^t \exp\left[-\g \int_0 ^\tau \|D^k v(\sigma )\|_{L^2}^{\frac{5}{2k}}d\sigma\right] d\tau,
$$
we find that (\ref{diff2}) can be rewritten as a differential  inequality,
\bb\label{diffeq2}
y'(t)\geq -\left(\g-C_0  \|v_0\|_{L^2}^{1-\frac{5}{2k}}\right)
\|D^k v_0 \|_{L^2}^{\frac{5}{2k}} \, y(t)^{M},
\ee
where $M$ is the same constant defined in (\ref{diffeq2}).
 The differential inequality  (\ref{diffeq2}) is solved as
\bb\label{y2}
y(t)\geq \frac{1}{\left(1+ C_0 \|v_0\|_{L^2}^{1-\frac{5}{2k}}\|D^k v_0 \|_{L^2}^{\frac{5}{2k}} t\right)^{\frac{\g}{C_0\|v_0\|_{L^2}^{1-\frac{5}{2k}}}-1}},
\ee
which proves  (\ref{thm2a}). $\square$\\
\ \\
\noindent{\bf Proof  of Theorem 1.2 } Let $T$ be the maximal time of existence of a classical solution $v$ of $(NS)_\nu$ in $L^p (\Bbb R^3)$, and $v\in C([0, T); L^p (\Bbb R^3 ))$.
For a solution $v(x,t)$ and the associated pressure function $\textbf{p}(x,t)$, we
define a functional transform from
$(v,\textbf{p})$ to $(V, P)$ defined by the formula,
\bq\label{7.4}
 v(x,t)&=&\exp\left[\frac{\g}{2}\int_0 ^t\| v(\tau)\npn d\tau \right] \ V\left(y, s \right),\\
 \label{7.5}
 \textbf{p}(x,t)&=&\exp\left[ \g\int_0 ^t\|v(\tau)\npn d\tau \right] \ P\left(y, s\right)
 \eq
 with
\bq \label{7.6}
y&=&\exp\left[\pm \frac{\g}{2}\int_0 ^t\|v(\tau)\npn d\tau \right]x,  \\
\label{7.6a}
 s&=&\int_0 ^t\exp\left[\g\int_0 ^\tau\| v(\sigma)\npn d\sigma \right]d\tau .
\eq
Here our choice of similarity transform makes the scaling dimension of the $\|v\|_{L^3}$ become zero, which is the natural choice for the (viscous) Navier-Stokes equations.  As a consequence we  have the following invariant of the transform,
$$ \|v(t)\|_{L^3}=\|V(s)\|_{L^3}.
$$
We also note the following integral invariant of the transform,
$$
\int_0 ^t \|  v(\tau)\npn  d\tau =\int_0 ^s \|  V(\sigma )\npn d\sigma \qquad  3<p\leq \infty.
$$
Substituting $(v,\textbf{p})$ in (\ref{7.4})-(\ref{7.6a}) into  $(NS)_\nu$, we  obtain an equivalent system of equations:
$$
 (NS)_*\left\{ \aligned
 &-\frac{\g }{2} \|V(s)\npn \left[  V + (y \cdot \nabla)V \right]=V_s+(V\cdot \nabla )V +\nabla P \\
& \qquad \hspace{1.in} -\nu \Delta V,\\
 & \mathrm{div}\, V=0,\\
&V(y,0)=V_0(y)=v_0 (y).
\endaligned \right.
$$
 Similarly to the above proof we observe that $V\in C([0, S ); L^p(\Bbb R^3))$, where
 $$S:=\int_0 ^T\exp\left[\g\int_0 ^\tau\|  v(\sigma)\npn d\sigma \right]d\tau
 $$
 is the maximal time of existence of the classical solution in $L^p (\Bbb R^3 )$ for the system $(NS)_\nu$.
 Operating div $(\cdot )$ on the first equations of  $(NS)_*$, we find
 $ -\Delta P=\mathrm{div}\, \mathrm{div}\, v\otimes v$, which implies the pressure-velocity relation,
 \bb\label{cz}
 P=\sum_{j,k=1}^3(-\Delta )^{\frac12}\partial_j (-\Delta )^{\frac12}\partial_k V_j V_k =\sum_{j,k=1}^3 R_j R_k V_jV_k,
 \ee
 which is well-known in the case of the original  Navier-Stokes equations $(NS)_\nu$, where $R_j, j=1,2,3$, is the Riesz transform in $\Bbb R^3$.
 Taking $L^2 (\Bbb R^3)$ inner product of the first equations of $(NS)_\nu $ with $V |V|^{p-2}$, and integrating by part, we find that
\bq\label{est7}
\lefteqn{\frac1p\frac{d}{ds} \|V(s)\|_{L^p}^p+\npg \| V\|_{L^p} ^{p+\np}+\frac{2\nu }{p} \|\nabla ( |V|^{\frac{p}{2}}) \|_{L^2}^2 }\hspace{.0 in}\n \\
&&= -\int_{\Bbb R^3} |V|^{p-2} (V\cdot \nabla ) P dy= \int_{\Bbb R^3} P(V\cdot \nabla )(|V|^{p-2})  dy\n \\
&& =\int_{\Bbb R^3} P(V\cdot \nabla ) \left(|V|^{\frac{p}{2}}\right)^{2-\frac{4}{p}} dy=\left(2-\frac{4}{p}\right) \int_{\Bbb R^3} P |V|^{\frac{p}{2}-2} (V\cdot \nabla ) (|V|^{\frac{p}{2}}) dy\n \\
&&\leq \left(2-\frac{4}{p}\right)\int_{\Bbb R^3} |P| \,|V|^{\frac{p}{2}-1} \left|\nabla  (|V|^{\frac{p}{2}})\right|dy
 \leq \left(2-\frac{4}{p}\right)\|P\|_{L^p} \|V\|_{L^p} ^{\frac{p}{2}-1}\| \nabla  (|V|^{\frac{p}{2}})\|_{L^2}\n \\
&& \leq C \|V\|_{L^p}^{\frac{p-1}{2}}\| \nabla  (|V|^{\frac{p}{2}})\|_{L^2}^{1+\frac{3}{p}}
\leq  \frac{2\nu }{p} \|\nabla ( |V|^{\frac{p}{2}}) \|_{L^2}^2 +C_0\| V\|_{L^p} ^{p+\np}
\eq
 for a constant $C_0=C_0( p, \nu)$, where we used the following estimate of the pressure,
 \bq\label{pres}
 \|P\|_{L^p} &\leq& C_p \|V\|_{L^{2p}} ^2 \leq C_p \|V\|_{L^p} ^{\frac12} \|V\|_{L^{3p}}^{\frac32}
 = C_p \|V\|_{L^p} ^{\frac12}\|\,|V|^{\frac{p}{2}}\, \|_{L^6} ^{\frac{3}{p}} \n \\
 &\leq & C_p \|V\|_{L^p} ^{\frac12}\|\nabla (|V|^{\frac{p}{2}}) \|_{L^2} ^{\frac{3}{p}}.
 \eq
The first estimate of (\ref{pres}) is due to the Calderon-Zygmund inequality applied to (\ref{cz}), while the last one follows by  applying the Sobolev imbedding $\dot{H}^1 (\Bbb R^3) \hookrightarrow L^6 (\Bbb R^3)$.
We also note that to get the first line of (\ref{est7}) we used the computations,
 \bqn
 &&\int_{\Bbb R^3} |V |^{p-2}[(y\cdot \nabla ) V ] \cdot V \,dy =\frac{1}{p}\int_{\Bbb R^3} (y\cdot \nabla )|V|^p dy\\
 &&\qquad=-\frac{1}{p}\int_{\Bbb R^3}[\mathrm{ div}\, y ] |V|^p dy =-\frac{3}{p}\|V\|_{L^p}^p.
 \eqn
 Absorbing the term $ \frac{2\nu }{p} \|\nabla ( |V|^{\frac{p}{2}}) \|_{L^2}^2$ to the left hand side
 in (\ref{est7}), we have the differential inequality
$$
\frac{d}{ds} \|V\|_{L^p}\leq -\left[\npg -C_0\right] \|V\|_{L^p} ^{1+\np},
$$
which can be solved to provide us with
\bb\label{sol1}
\| V(s)\|_{L^p}\leq \frac{\|V_0\|_{L^p}}{\left[1+\left(\g-C_0\right)
\| V_0 \npn s\right]^{\npi} }
\ee
for all $s\in [0, S)$.
Transforming back to the original velocity $v$, using the relations (\ref{7.4})-(\ref{7.6a}),  we  obtain (\ref{thm7}).
In order to derive (\ref{thm7a}) we rewrite (\ref{thm7}) in the integrable form,
\bq\label{log7}
&&\| v(t)\npn \leq \frac{\| v_0\npn \exp \left[ \g \int_0 ^t \| v (\tau)\npn d\tau\right]}{\left\{1+\left(\g-C_0 \right)
\| v_0 \npn \int_0 ^t \exp\left[\g \int_0 ^\tau \| v(\sigma )\npn d\sigma\right] d\tau \right\}}\n \\
&&=\left(\g -C_0 \right)^{-1} \times \n \\
&&\times \frac{d}{dt}\log \left\{1+\left(\g-C_0 \right)\| v_0 \npn \int_0 ^t \exp\left[\g \int_0 ^\tau \| v(\sigma )\npn d\sigma\right] d\tau \right\}.\n \\
\eq
Hence, integrating  (\ref{log7}) over $[0, t]$, we obtain
\bq\label{diff7}
&&\int_0 ^t\|v(\tau)\npn d\tau\leq \left(\g-C_0 \right)^{-1}\times\n \\
&&\times \log \left\{1+\left(\g-C_0 \right)
\| v_0 \npn \int_0 ^t \exp\left[\g \int_0 ^\tau \| v(\sigma )\npn d\sigma\right] d\tau \right\}.\n \\
\eq
Now, setting
$$y(t):= 1+\left(\g-C_0\right)
\|v_0 \npn \int_0 ^t \exp\left[\g \int_0 ^\tau \| v(\sigma )\npn d\sigma\right] d\tau,
$$
we find that (\ref{diff7}) can be rewritten as a differential  inequality,
\bb\label{diffeq7}
y'(t)\leq \left(\g-C_0 \right)
\| v_0 \npn \, y(t)^{\frac{\g}{\g-C_0}},
\ee
which can be solved as
\bb\label{y7}
y(t)\leq \frac{1}{\left(1- C_0 \|v_0 \npn t\right)^{\frac{\g}{C_0}-1}},
\ee
which provides us with (\ref{thm7a}). $\square$\\
\ \\
\noindent{\bf Proof  of Theorem 1.3 }
 Let $T$ be the maximal time of existence of a classical solution $\t $ of $(QG)_\nu$ in $W^{m,p}(\Bbb R^2)$, and $\t \in C([0, T); W^{m,p}(\Bbb R^2))$.
This time we introduce a functional transform from
$(\t, v) $ to $(\Theta, V)$ defined by the formula,
\bq\label{3.4}
 \t (x,t)&=&\exp\left[\pm\frac{2\g}{p+2}\int_0 ^t\| D^k \t (\tau)\|_{L^p} ^{\frac{p+2}{kp}} d\tau \right]  \T\left(y, s \right),\\
 \label{3.5}
 v(x,t)&=&\exp\left[\pm\frac{2\g}{p+2}\int_0 ^t\| D^k \t (\tau)\|_{L^p} ^{\frac{p+2}{kp}} d\tau \right]  V \left(y, s \right)
 \eq
 with
\bq \label{3.6}
y&=&\exp\left[\pm \frac{p\g}{p+2}\int_0 ^t\| D^k \t (\tau)\|_{L^p} ^{\frac{p+2}{kp}} d\tau  \right]x,  \\
\label{3.6a}
 s&=&\int_0 ^t\exp\left[\pm \g\int_0 ^\tau\| D^k \t (\sigma)\|_{L^p} ^{\frac{p+2}{kp}} d\sigma \right]d\tau ,
\eq
respectively for $(\pm)$.
Here we notice that our choice of similarity transform makes the scaling dimension of $\|\theta(t)\|_{L^p}$ become zero, and we have the invariants of the transform,
$$
\|\theta(t)\|_{L^p}=\|\Theta (s)\|_{L^p}\qquad   0<p\leq \infty ,
$$
and
$$
\int_0 ^t\| D^k \t (\tau)\pn d\tau= \int_0 ^s\| D^k \T (\sigma)\pn d\sigma .
$$
Substituting $(v,p)$ in (\ref{3.4})-(\ref{3.6a}) into  $(QG)_\k$, we  obtain an equivalent system of equations:
$$
 (QG)_\k ^\pm \left\{ \aligned
 &\mp\frac{\g }{p+2} \|D^k \T (s)\|_{L^p} ^{\frac{p+2}{kp}} \left[ 2 \T + p (y \cdot \nabla)\T \right]=\T_s+(V\cdot \nabla )\T \\
& \qquad \hspace{.3in} -\k \Lambda ^\a\T \exp \left[ \mp\left(1-\frac{p\a}{p+2} \right)\g\int_0 ^s \| D^k \T (\sigma)\|_{L^p} ^{\frac{p+2}{kp}}d\sigma \right],\\
 & V= \nabla ^\bot (-\Delta )^{-\frac12}\T  ,\\
&\T(y,0)=\T_0(y)=\t_0 (y),
\endaligned \right.
$$
  where $(QS)_\k ^+$ means that we have chosen (+) sign in (\ref{3.4})-(\ref{3.6a}), and this corresponds to $(-)$  sign in the first equations of $(QG)_\k ^\pm $.  Similarly for $(QG)_\k ^-$.
Similarly to the proof of Theorem 1.1 we observe that $\T \in C([0, S_\pm ); W^{m,p}(\Bbb R^2))$, where
 $$S_\pm:=\int_0 ^T\exp\left[\pm\g\int_0 ^\tau\| D^k \t (\sigma)\pn d\sigma \right]d\tau
 $$
 is the maximal time of existence of the classical solution in $W^{m,p}(\Bbb R^2)$ for the system $(QG)_\k ^\pm$ respectively.\\

 \noindent{\it \underline{Proof of (i)}: } We choose $(+)$ sign in (\ref{3.4})-(\ref{3.6a}), and
 work with $(QG)_\k ^+$, where $\k \geq 0$.
Taking $L^2(\Bbb R^2)$ inner product of the first equations of $(QG)_\k ^+$ by $\T |\T |^{p-2}$, and integrating by part, we find that
\bqn
 \lefteqn{\frac1p\frac{d}{ds} \|\T (s)\|_{L^p} ^p=
 -k\int_{\Bbb R^2} \T |\T|^{p-2} \la ^\a \T \,dy}\hspace{.1in}\n \\
 &\leq &-\frac{2\k }{p} \int_{\Bbb R^2} \left|\la^{\frac{\a}{2}}\left(|\T
 |^{\frac{p}{2}}\right) \right|^2 dy\, \exp \left[ -\frac{\g}{5}\int_0 ^s \|D^k \T(\sigma)\pn d\sigma \right]\leq 0,
 \eqn
 where we used (\ref{frac}) for the viscosity term.
 Thus, we have the $L^p$ bound of $\T$.
\bb
 \|\T (s)\|_{L^p}\leq \|\T_0\|_{L^p},
 \ee
Next, operating $D^k $ on the first equations of $(QG)_\k ^+ $ and then taking $L^2 (\Bbb R^2)$ inner product of it by $D^k \T |D^k \T |^{p-2}$, and integrating by part,
we estimate
\bq\label{est3}
\lefteqn{\frac1p \frac{d}{ds} \|D^k \T \|_{L^p} ^p +\frac{kp\g}{p+2} \|D^k \T\|_{L^p} ^{2+\frac{p+2}{kp}} }\hspace{.2 in}\n \\
 &&+\frac{2\k}{p} \int_{\Bbb R^2} |\la ^{\frac{\a}{2}} (D^k \T)^2 |^{\frac{p}{2}}dy \exp \left[- \left(1-\frac{p\a}{p+2} \right)\int_0 ^s \|D^k \T (\sigma)\|_{L^2}^{\frac{p+2}{kp}}d\sigma \right]\n \\
&&\leq - \int_{\Bbb R^2}\left[D^k (V\cdot \nabla )\T- (V\cdot \nabla )D^k \T \right]D^k\T|D^k \T|^{p-2} dy \n \\
&&\leq \|D^k (V\cdot \nabla ) \T- (V\cdot \nabla )D^k \T\|_{L^p} \|D^k \T\|^{p-1}_{L^p} \n \\
&& \leq C (\|\nabla V\|_{L^\infty} +\|\nabla \T \|_{L^\infty}) ( \|D^k V\|_{L^p} +\|D^k \T \|_{L^p})\|D^k \T\|_{L^p} ^{p-1}\n \\
&&\leq C  (\|V\|_{L^p}^{1-\frac{p+2}{kp}}\|D^k V\|_{L^p} ^{\frac{p+2}{kp}}+ \|\T \|_{L^p}^{1-\frac{p+2}{kp}}\|D^k \T \|_{L^p} ^{\frac{p+2}{kp}}) \|D^k \T \|_{L^p}^p\n \\
&&\leq \frac{kp C_0}{p+2} \|\T \|_{L^p}^{1-\frac{p+2}{kp}} \|D^k \T\|_{L^p} ^{2+\frac{p+2}{kp}} \leq \frac{kp C_0}{p+2} \|\T_0 \|_{L^p}^{1-\frac{p+2}{kp}} \|D^k \T \|_{L^p} ^{2+\frac{p+2}{kp}}
\eq
for an absolute constant $C_0=C_0 (k, p)$.
In (\ref{est3}) we used the computation,
\bqn
&&\int_{\Bbb R^2} D^k [(y\cdot \nabla )\T ]D^k \T |D^k \T|^{p-2} dy = \frac1p \int_{\Bbb R^2} (y\cdot \nabla )|D^k \T|^p dy + k\|D^k \T\|_{L^p}^p\n \\
&&=-\frac2p \|D^k \T\|_{L^p}^p + k\|D^k \T\|_{L^p}^p=\left(k-\frac2p\right) \|D^k \T\|_{L^p}^p
\eqn
the commutator estimate (\ref{com}) and the Gagliardo-Nirenberg inequality  (\ref{gag}), and also
the Calderon-Zygmund type of inequality,
$$
\|V\|_{\dot{W}^{k,p}} \leq C \|\T \|_{\dot{W}^{k,p}}, \qquad 1<p<\infty, \quad k\in \Bbb N \cup \{0\}.
$$
Hence, from (\ref{est3}), ignoring the viscosity term, we have the differential inequality
$$
\frac{d}{ds} \|D^k \T \|_{L^p}\leq -\frac{kp}{p+2}\left(\g -C_0\|\T_0\|_{L^p}^{1-\frac{p+2}{kp}}\right) \|D^k \T\|_{L^p} ^{1+\frac{p+2}{kp}},
$$
which can be solved to provide us with
\bb\label{sol3}
\|D^k \T (s)\|_{L^p}\leq \frac{\|D^k \T_0\|_{L^p}}{\left[1+\left(\g-C_0 \|\T_0\|_{L^p}^{1-\frac{p+2}{kp}}\right)
\|D^k \T_0 \|_{L^p}^{\frac{p+2}{kp}} s\right]^{\frac{kp}{p+2}} }
\ee
for all $s\in [0, S_+)$.
Transforming back to the original velocity $v$, using the relations (\ref{3.4})-(\ref{3.6a}) with $(+)$ sign,  we derive (\ref{thm3}).
We now derive (\ref{thm3a}). For this we note that  (\ref{thm3}) can be written in the  integrable form,
\bqn\label{log3}
\lefteqn{\|D^k v(t)\|_{L^p}^{\frac{p+2}{kp}} \leq \frac{\|D^k \t_0\|_{L^p}^{\frac{p+2}{kp}} \exp \left[ \g \int_0 ^t \|D^k \t (\tau)\|_{L^p} ^{\frac{p+2}{kp}} d\tau\right]}{\left\{1+\left(\g-C_0  \|\t_0\|_{L^p}^{1-\frac{(p+2)}{kp}}\right)
\|D^k \t_0 \|_{L^p}^{\frac{p+2}{kp}}\int_0 ^t \exp\left[\g \int_0 ^\tau \|D^k \t(\sigma )\|_{L^p}^{\frac{p+2}{kp}}d\sigma\right] d\tau \right\}}}\hspace{.0in}\n \\
&&=\left(\g -C_0  \|\t_0\|_{L^p}^{1-\frac{p+2}{kp}}\right)^{-1} \times \n \\
&&\times \frac{d}{dt}\log \left\{1+\left(\g-C_0 \|\t_0\|_{L^p}^{1-\frac{p+2}{kp}}\right)\|D^k \t_0 \|_{L^p}^{\frac{p+2}{kp}}\int_0 ^t \exp\left[\g \int_0 ^\tau \|D^k \t(\sigma )\|_{L^p}^{\frac{p+2}{kp}}d\sigma\right] d\tau \right\}.\n \\
\eqn
Integrating this  over $[0, t]$, we obtain
\bq\label{diff3}
\lefteqn{\int_0 ^t\|D^k \t(\tau)\|_{L^p}^{\frac{p+2}{kp}}d\tau\leq \left(\g-C_0 \|\t_0\|_{L^p}^{1-\frac{p+2}{kp}}\right)^{-1}\times}\hspace{.0in}\n \\
&&\times \log \left\{1+\left(\g-C_0   \|\t_0\|_{L^p}^{1-\frac{(p+2)}{kp}}\right)
\|D^k \t_0 \|_{L^p}^{\frac{p+2}{kp}}\int_0 ^t \exp\left[\g \int_0 ^\tau \|D^k \t(\sigma )\|_{L^p}^{\frac{p+2}{kp}}d\sigma\right] d\tau \right\}.\n \\
\eq
Setting
$$y(t):= 1+\left(\g-C_0  \|\t_0\|_{L^p}^{1-\frac{p+2}{kp}}\right)
\|D^k \t_0 \|_{L^p}^{\frac{p+2}{kp}}\int_0 ^t \exp\left[\g \int_0 ^\tau \|D^k \t(\sigma )\|_{L^p}^{\frac{p+2}{kp}}d\sigma\right] d\tau,
$$
we find that (\ref{diff3}) can be rewritten as a differential  inequality,
\bb\label{diffeq3}
y'(t)\leq \left(\g-C_0  \|\t_0\|_{L^p}^{1-\frac{p+2}{kp}}\right)
\|D^k \t_0 \|_{L^p}^{\frac{p+2}{kp}} \, y(t)^{M},
\ee
where we set
\bb\label{m3}
M:=\frac{\g}{\g-C_0  \|\t_0\|_{L^p}^{1-\frac{p+2}{kp}}}.
\ee
 The differential inequality  (\ref{diffeq3}) is solved as
\bb\label{y3}
y(t)\leq \frac{1}{\left(1- C_0 \|\t_0\|_{L^p}^{1-\frac{p+2}{kp}}\|D^k \t_0 \|_{L^p}^{\frac{p+2}{kp}} t\right)^{\frac{\g}{C_0\|\t_0\|_{L^p}^{1-\frac{p+2}{kp}}}-1}},
\ee
which provides us with (\ref{thm3a}).\\
\ \\
\noindent{\it \underline{Proof of (ii)}: }  We choose $(-)$ sign in (\ref{3.4})-(\ref{3.6a}), and
 work with $(QG)_0^-$.
Taking $L^2(\Bbb R^2)$ inner product of the first equations of $(QG)_0 ^+$ by $\T |\T |^{p-2}$, and integrating by part, we find  first that
 $$
 \frac1p\frac{d}{ds} \|\T (s)\|_{L^p} ^p=0,
 $$
  which implies the $L^p$ conservation of  $\T$.
\bb
 \|\T (s)\|_{L^p}=\|\T_0\|_{L^p},
 \ee
Next, operating $D^k $ on the first equations of $(QG)_0 ^- $ and then taking $L^2 (\Bbb R^2)$ inner product of it by $D^k \T |D^k \T |^{p-2}$, and integrating by part,
we estimate below
\bq\label{est4}
\lefteqn{\frac1p \frac{d}{ds} \|D^k \T \|_{L^p} ^p -\frac{kp\g}{p+2} \|D^k \T\|_{L^p} ^{2+\frac{p+2}{kp}} }\hspace{.2 in}\n \\
&&=- \int_{\Bbb R^2}\left[D^k (V\cdot \nabla )\T- (V\cdot \nabla )D^k \T \right]D^k\T|D^k \T|^{p-2} dy \n \\
&&\geq -\|D^k (V\cdot \nabla ) \T- (V\cdot \nabla )D^k \T\|_{L^p} \|D^k \T\|^{p-1}_{L^p} \n \\
&& \geq -C (\|\nabla V\|_{L^\infty} +\|\nabla \T \|_{L^\infty}) ( \|D^k V\|_{L^p} +\|D^k \T \|_{L^p})\|D^k \T\|_{L^p} ^{p-1}\n \\
&&\geq- C  (\|V\|_{L^p}^{1-\frac{p+2}{kp}}\|D^k V\|_{L^p} ^{\frac{p+2}{kp}}+ \|\T \|_{L^p}^{1-\frac{p+2}{kp}}\|D^k \T \|_{L^p} ^{\frac{p+2}{kp}}) \|D^k \T\|_{L^p}^p\n \\
&&\geq -\frac{kpC_0}{p+2}\|\T \|_{L^p}^{1-\frac{p+2}{kp}} \|D^k \T\|_{L^p} ^{2+\frac{p+2}{kp}} \geq - \frac{kpC_0}{p+2}\|\T_0 \|_{L^p}^{1-\frac{p+2}{kp}} \|D^k \T\|_{L^p} ^{2+\frac{p+2}{kp}}\n \\
\eq
for the same  absolute constant $C_0=C_0 (k, p)$ as in the proof of (i).
Hence, from (\ref{est4}),  we have the differential inequality
$$
\frac{d}{ds} \|D^k \T \|_{L^p}\geq  \frac{kp}{p+2}\left(\g -C_0\|\T_0\|_{L^p}^{1-\frac{p+2}{kp}}\right) \|D^k \T\|_{L^p} ^{1+\frac{p+2}{kp}},
$$
which can be solved to provide us with
\bb\label{sol4}
\|D^k \T (s)\|_{L^p}\geq \frac{\|D^k \T_0\|_{L^p}}{\left[1-\left(\g-C_0 \|\T_0\|_{L^p}^{1-\frac{p+2}{kp}}\right)
\|D^k \T_0 \|_{L^p}^{\frac{p+2}{kp}} s\right]^{\frac{kp}{p+2}} }
\ee
for all $s\in [0, S_-)$.
Transforming back to the original velocity $v$, using the relations (\ref{3.4})-(\ref{3.6a}) with $(-)$ sign,  we obtain (\ref{thm4}).
In order to prove (\ref{thm4a}) we note that (\ref{thm4})  can be written in the  integrable form,
\bqn\label{log4}
\lefteqn{\|D^k v(t)\|_{L^p}^{\frac{p+2}{kp}} \geq }\hspace{.0in}\n \\
&&\geq \frac{\|D^k \t_0\|_{L^p}^{\frac{p+2}{kp}} \exp \left[ -\g \int_0 ^t \|D^k \t (\tau)\|_{L^p} ^{\frac{p+2}{kp}} d\tau\right]}{\left\{1-\left(\g-C_0  \|\t_0\|_{L^p}^{1-\frac{(p+2)}{kp}}\right)
\|D^k \t_0 \|_{L^p}^{\frac{p+2}{kp}}\int_0 ^t \exp\left[-\g \int_0 ^\tau \|D^k \t(\sigma )\|_{L^p}^{\frac{p+2}{kp}}d\sigma\right] d\tau \right\}}\\
&&=-\left(\g -C_0  \|\t_0\|_{L^p}^{1-\frac{p+2}{kp}}\right)^{-1} \times \n \\
&&\times \frac{d}{dt}\log \left\{1-\left(\g-C_0 \|\t_0\|_{L^p}^{1-\frac{p+2}{kp}}\right)\|D^k \t_0 \|_{L^p}^{\frac{p+2}{kp}}\int_0 ^t \exp\left[\g \int_0 ^\tau \|D^k \t(\sigma )\|_{L^p}^{\frac{p+2}{kp}}d\sigma\right] d\tau \right\}.\n \\
\eqn
Integrating this  over $[0, t]$, we obtain
\bq\label{diff4}
\lefteqn{\int_0 ^t\|D^k \t(\tau)\|_{L^p}^{\frac{p+2}{kp}}d\tau\geq - \left(\g-C_0 \|\t_0\|_{L^p}^{1-\frac{p+2}{kp}}\right)^{-1}\times}\hspace{.0in}\n \\
&&\times \log \left\{1-\left(\g-C_0   \|\t_0\|_{L^p}^{1-\frac{(p+2)}{kp}}\right)
\|D^k \t_0 \|_{L^p}^{\frac{p+2}{kp}}\int_0 ^t \exp\left[-\g \int_0 ^\tau \|D^k \t(\sigma )\|_{L^p}^{\frac{p+2}{kp}}d\sigma\right] d\tau \right\}.\n \\
\eq
Setting
$$y(t):= 1-\left(\g-C_0  \|\t_0\|_{L^p}^{1-\frac{p+2}{kp}}\right)
\|D^k \t_0 \|_{L^p}^{\frac{p+2}{kp}}\int_0 ^t \exp\left[\g \int_0 ^\tau \|D^k \t(\sigma )\|_{L^p}^{\frac{p+2}{kp}}d\sigma\right] d\tau,
$$
we find that (\ref{diff4}) can be rewritten as a differential  inequality,
\bb\label{diffeq4}
y'(t)\geq -\left(\g-C_0  \|\t_0\|_{L^p}^{1-\frac{p+2}{kp}}\right)
\|D^k \t_0 \|_{L^p}^{\frac{p+2}{kp}} \, y(t)^{M},
\ee
where
$$
M=\frac{\g}{\g-C_0  \|\t_0\|_{L^p}^{1-\frac{p+2}{kp}}}.
$$
 The differential inequality  (\ref{diffeq4}) is solved as
\bb\label{y4}
y(t)\geq \frac{1}{\left(1+ C_0 \|\t_0\|_{L^p}^{1-\frac{p+2}{kp}}\|D^k \t_0 \|_{L^p}^{\frac{p+2}{kp}} t\right)^{\frac{\g}{C_0\|\t_0\|_{L^p}^{1-\frac{p+2}{kp}}}-1}},
\ee
which provides us with (\ref{thm4a}). $\square$\\
\ \\
\noindent{\bf Proof of Theorem 1.4 } We transform from $(\t, v)$ to $(\T, V)$ according to the formula
\bq\label{4.1}
 \t(x,t)&=&\exp\left[\frac{\pm\g\laa}{\laa+1}\int_0 ^t\|\nt (\tau)\bn d\tau \right] \ \T\left(y, s \right),\\
 \label{4.2}
 v(x,t)&=&\exp\left[\frac{\pm \g\laa}{\laa+1}\int_0 ^t\|\nt (\tau)\bn d\tau \right] \ V \left(y, s\right)
 \eq
 with
\bq \label{4.3}
y&=&\exp\left[\frac{\pm \g}{\laa+1}\int_0 ^t\|\nt(\tau)\bn d\tau \right]x, \n \\
 s&=&\int_0 ^t\exp\left[\pm \g\int_0 ^\tau\|\nt (\sigma)\bn d\sigma \right]d\tau
\eq
respectively for the signs $\pm$. In (\ref{4.1})-(\ref{4.3}) both $\g>0$ and $\laa >-1$ are free parameters. We note the following integral invariant,
$$
\int_0 ^t\|\nt (\tau)\bn d\tau=\int_0 ^s\|\nO (\sigma)\bn d\sigma
$$
for all $\laa >-1$.
Substituting (\ref{4.1})-(\ref{4.3}) into the $(QG)_0$, we find that
 $$
 (QG_\pm) \left\{ \aligned
  & \mp \g \|\nabla \T (s)\bn\left[ \frac{\laa}{\laa+1} \T  +\frac{1}{\laa+1}(y \cdot \nabla)\T \right]=\T_s+(V\cdot \nabla )\T ,\\
 & V=\nabla ^\bot (-\Delta )^{-\frac12} \T ,\\
 & \T(y,0)=\T_0(y)=\t_0 (y)
  \endaligned \right.
 $$
 respectively for $\pm$.
Below we denote $( \T^\pm, V^\pm )$ for the solutions of $(QG_\pm)$ respectively.
We will first derive the following estimates for the system $(QG_\pm)$.
 \bq\label{4.4}
\| \nO^+(s)\|_{L^\infty} &\leq& \frac{\|\nO_0\|_{L^\infty}}{1+(\g-1 )s\|\nO _0 \|_{L^\infty}},\\
\label{4.5}
\| \nO ^-(s)\|_{L^\infty}&\geq &\frac{\|\nO_0\|_{L^\infty}}{1-(\g-1 )s\|\nO_0 \|_{L^\infty}},
\eq
as long as $\T ^\pm(s)\in B^{1}_{\infty, 1} (\Bbb R^3 )$.
Taking operation of $\nabla ^\bot$ on  the first equation of $(QG_\pm)$, we have
\bb\label{4.6}
\mp\g \|\nO\bn \left[\nOo -\frac{1}{\laa +1}(y\cdot \nabla )\nOo\right]=
\nOo_s +(V\cdot \nabla )\nOo-(\nOo\cdot \nabla )V.
\ee
Multiplying $\Xi=\nOo/|\nOo|$ on the both sides of (\ref{4.6}), we deduce
\bq\label{4.7}
\lefteqn{|\nO|_s +(V\cdot \nabla )|\nO|\mp\frac{ \|\nO \bn}{\laa +1}(y\cdot \nabla )|\nO| }\hspace{.2in}\n \\
&&=(\Xi\cdot \nabla V\cdot \Xi \mp C_0 \|\nO \bn |\nO|)\n \\
&&\qquad\mp(\g -C_0)\|\nO\bn |\nO|\n \\
&&\left\{ \aligned &
 \leq -(\g -C_0)\|\nO \bn |\nO|\quad \mbox{for}\,(QG_+),\\
 &\geq (\g -C_0)\|\nO \bn |\nO|\quad \mbox{for}\, (QG_-),
 \endaligned \right.
\eq
where we used the estimates
$$|\Xi\cdot \nabla V\cdot \Xi |\leq |\nabla V|\leq \|\nabla V\|_{L^\infty} \leq \|\nV\bn
\leq C_0 \|\nO \bn  $$
for an absolute constant $C_0$, the last step of which follows by the Calderon-Zygmund type of inequality
on $\dot{B}^{0}_{\infty,1} (\Bbb R^2)$.
Given smooth solution pair $(V, \T )$ of the system $(QG_\pm)$, we introduce the particle trajectories $\{ Y(a,s)= Y_\pm (a,s)\}$ defined by
$$
\frac{\partial Y(a,s)}{\partial s}=V(Y(a,s),s)\mp \frac{\|\nO  (s)\bn}{\laa+1} Y(a,s) \quad ;\quad
Y(a,0)=a.
$$
Using  the inequalities
$$ \|\nO \bn \geq \|\nO \|_{L^\infty}\geq  |\nO (y,s)|\qquad \forall y\in \Bbb R^3,
$$
we can further estimate from (\ref{4.7})
\bb\label{4.8}
\frac{\partial}{\partial s} |\nO(Y(a,s),s)| \left\{ \aligned &
 \leq -(\g -C_0)|\nO(Y(a,s),s)|^2 \quad \mbox{for}\, (QG_+),\\
 &\geq (\g -C_0) |\nO(Y(a,s),s)|^2\quad \mbox{for}\, (QG_-) .
 \endaligned \right.
\ee
Solving these differential inequalities (\ref{4.7}) along the particle trajectories, we obtain that
\bb\label{4.9}
 |\nO (Y(a,s) ,s) | \left\{ \aligned &
 \leq \frac{|\nO_0 (a)|}{1 +(\g -C_0)s|\nO_0(a )|}\quad \mbox{for}\,(QG_+)\\
 &\geq \frac{|\nO_0 (a)|}{1 -(\g -C_0)s|\nO_0(a )|}\quad \mbox{for}\, (QG_-).
 \endaligned \right.
\ee
Writing the first inequality of (\ref{4.9}) as
$$
 |\nO^+ (Y(a,s) ,s) | \leq \frac{1}{\frac{1}{|\nO_0 (a)|} +(\g -C_0 )s}
 \leq \frac{1}{\frac{1}{\|\nO_0 \|_{L^\infty}} +(\g -C_0)s},
 $$
 and then taking supremum over $a\in \Bbb R^2$, which is equivalent to taking supremum
 over $Y(a,s)\in \Bbb R^2$ due to the fact that the mapping $a\mapsto Y(a,s)$ is a deffeomorphism on $\Bbb R^2$ as long as $V\in C([0, S); \dot{B}^1_{\infty, 1} (\Bbb R^2))$, we obtain (\ref{4.4}).
 In order to derive (\ref{4.5})  from the second inequality of (\ref{4.9}), we first write
 $$
 \|\nO^-( s)\|_{L^\infty}\geq |\nO (Y(a,s),s)|\geq  \frac{1}{\frac{1}{|\nO_0 (a) |} -(\g -C_0)s},
 $$
 and than take supremum over $a\in \Bbb R^2$.
Finally, in order  to obtain  (\ref{thm5})-(\ref{thm6}), we just change variables from (\ref{4.4})-(\ref{4.5})  back to the original physical ones, using the fact
\bqn
\nO^+ (y,s)&=&\exp\left[ -\g \int_0 ^t \|\nt(\tau )\bn d\tau\right]\n (x,t),\\
s&=&\int_0 ^t \exp\left[ \g \int_0 ^\tau \|\nt(\sigma )\bn d\sigma\right]d\tau
\eqn
for (\ref{thm5}), while  in order to deduce (\ref{thm6}) from (\ref{4.5}) we substitute
\bqn
\nO^- (y,s)&=&\exp\left[ \g \int_0 ^t \|\nt(\tau)\bn d\tau \right]\o (x,t),\\
s&=&\int_0 ^t \exp\left[ -\g \int_0 ^\tau\|\nt(\sigma)\bn d\sigma\right]d\tau .
\eqn
In order to derive (\ref{thm6a}) we observe that (\ref{thm6}) can be written as
$$
 \|\nt (t)\|_{L^\infty} \geq \frac{-1}{(\g-C_0)} \frac{d}{dt} \left\{1-(\g-C_0) \|\nt_0 \|_{L^\infty}\int_0 ^t\exp\left[ -\g \int_0 ^\tau\|\nt (\sigma)\bn d\sigma \right]d\tau\right\},
 $$
 which, after integration over $[0,t]$, provides us with the estimates,
 \bq\label{last5}
&& \int_0 ^t \|\nt (\tau )\bn  d\tau\geq  \int_0 ^t \|\nt (\tau )\|_{L^\infty} d\tau \n \\
&& \geq \frac{-1}{(\g-C_0)}\log\left\{  1-(\g-C_0)\|\nt_0 \|_{L^\infty}\int_0 ^t\exp\left[ -\g \int_0 ^\tau\|\nt (\sigma)\bn d\sigma \right]d\tau\right\}\n \\
 \eq
 for all  $\g >C_0$.
 Setting
 $$ y(t):=1-(\g-C_0)\|\nt_0 \|_{L^\infty}\int_0 ^t\exp\left[ -\g \int_0 ^\tau\|\nt (\sigma)\bn d\sigma \right]d\tau,
 $$
 the inequality (\ref{last5}) can be written as another differential inequality,
 $$
 y'(t)\geq -(\g-C_0) \|\nt_0 \|_{L^\infty} y(t)^{\frac{\g}{\g-C_0}}.
 $$
 Solving this we obtain (\ref{thm6a}).
 
\end{document}